\documentclass[11pt,reqno]{amsart}
\usepackage{a4wide}
\usepackage[all]{xy}
\newcommand\bbN{{\mathbb N}}
\renewcommand\a{{\alpha}}
\newcommand\mfa{{\mathfrak a}}
\newcommand\mfb{{\mathfrak b}}
\newcommand\mfU{{\mathfrak U}}
\let\rar\rightarrow
\let\lar\longrightarrow
\let\disp\displaystyle

\newtheorem{stuff}{Stuff}[section]
\newtheorem{theorem}[stuff]{\bf Theorem}
\newtheorem{proposition}[stuff]{\bf Proposition}
\newtheorem{lemma}[stuff]{\bf Lemma}

\newtheorem{consequence}[stuff]{\bf Consequence}
\newenvironment{remark}{%
\vskip1ex\refstepcounter{stuff}\trivlist \itemindent 0pt
\item[\hskip\labelsep\bf Remark \thestuff.]%
\ignorespaces}{\endtrivlist\vskip1ex}%

\begin{document}
\title{Comparison of certain complexes of modules of generalized 
fractions and \v Cech complexes}

\author{R.~Tajarod, H.~Zakeri$^*$}
\address{Mosaheb Institute of Mathematics\newline  
\null\quad Teacher Training University,\newline 
\null\quad 599 Taleghani Ave., Tehran 15618, Iran}
\thanks{{\it Keywords}\,: \v Cech complex, local cohomology, modules of generalized fractions, filter regular sequence\newline 
\null\hspace{3ex}{\it E-mail}\,: tajarod\@@neda.net.ir, zakeri\@@saba.tmu.ac.ir\newline 
\null\quad$^*$Research supported in part by MIM Grant No. P83-117}
\subjclass[2000]{13D25, 13D45}

\begin{abstract}
We establish an explicit quasi-isomophism of complexes, 
which is homogeneous in graded situation, from a given 
\v Cech complex of modules of generalized fractions. 
\end{abstract}

\maketitle

\markboth{Comparison of generalized fractions and \v Cech complexes}
{R.~Tajarod, H.~Zakeri}

\section*{Introduction}

Throughout this note, $A$ is a commutative Noetherian ring with non-zero
identity and $M$ is a finitely generated $A$-module. 

One of the popular approaches to the construction of the local cohomology
modules uses cohomology of \v Cech complexes. In this note we establish explicit isomorphisms between the homology modules of a given \v Cech complex of $M$ and the corresponding homology modules of a certain complex of modules of generalized fractions. For a sequence $x_1,\ldots,x_n$ of elements of $A$, we shall use $C(x, M)^\bullet$ to denote the \v Cech complex of $M$ with respect to $x_1,\ldots, x_n$. The reader is referred, for example, to \cite[\S 5]{bs} for the construction of \v Cech complexes. As we mentioned above, this note is also concerned with a certain complex $C(\mathfrak U(x),M)$ of $A$-modules and  $A$-homomorphisms which involves modules of generalized fractions derived from $M$ and the sequence $x_1,\ldots,x_n$. The complex is described as follows: \newline for each $i\in\bbN$ (throughout we use $\bbN$ to denote the set of positive  integers) set 
$$
U(x)_i = \biggl\{
(x_1^{\a_1},... ,x_i^{\a_i}) 
\,\left|\, 
\begin{array}{r}
\text{there exists $j$ with $0\leq j\leq i$ 
such that $\a_1,... ,\a_j\in\bbN$,} \\ 
\text{and $\a_{j+1}=\ldots =\a_i = 0$ }
\end{array}
\right.
\biggr\},
$$
where $x_r$ is interpreted as $x_n$ whenever $r>n$. Then, for each $i\in\bbN$, $U(x)_i$ is a triangular subset \cite[2.1]{sz1} of $A^i$. We can \cite[page 420]{oc} form the complex 
$$ 
0\overset{e^{-1}}{\lar} 
M \overset{e^{0}}{\lar} U(x)_1^{-1}M 
\overset{e^{1}}{\lar} \ldots {\lar}
U(x)_i^{-i}M \overset{e^{i}}{\lar} U(x)_{i+1}^{-i-1}M{\lar}\ldots
$$ 
of $A$-modules and $A$-homomorphisms, which we denote by $C(\mathfrak U(x),M)$. Here $U(x)_i^{-i}M$ is the module of generalized fractions of $M$ with respect to the triangular subset $U(x)_i$ of $A^i$, and the homomorphisms $e^i$ (for $i\geq 0$) are given by the following formulas: 
\newline $e^0(b) = \frac{b}{(1)}$ for all $b\in M$ and, for each $i\in\bbN$, 
$$ 
e^i\biggl(
\frac{b}{(x_1^{\a_1},\ldots,x_i^{\a_i})}
\biggr)
= 
\frac{(-1)^ib}{(x_1^{\a_1},\ldots,x_i^{\a_i},1)}
$$
for all $b\in M$ and $(x_1^{\a_1},\ldots,x_i^{\a_i})\in U(X)_i$. 

It is well known that the homology modules of the complex $C(x, M)^\bullet$ are, up to $A$-isomorphisms, independent of the choice of sequence of generators for the ideal $\mathfrak a =\overset{n}{\underset{i=1}{\sum}}Ax_i$. 

Thus, to establish the main result of this note, without loss of generality, we may choose a better generating set for the ideal $\mathfrak a$. Indeed, with that choice of $x_1,\ldots,x_n$ (see proposition \ref{prop12}), we will establish a quasi-isomorphism from the complex $C(x, M)^\bullet$ to the complex $C(\mathfrak U(x), M)$, which is homogeneous in the graded situation. This result provides, perhaps, a simpler method to compute the homology modules of the the \v Cech complex $C(x, M)^\bullet$.

\section{The results}{\label{sct:results}}
 
Throughout this paper $\mfa$ is an ideal of $A$. For a submodule $N$ of $M$, we set 
$$ 
N:_M\mfa^\infty := \{m\in M \mid \mfa^rm\subseteq N\text{ for some }r\in\bbN\}.
$$
We recall the definition of a filter regular sequence. A sequence $x_1,\ldots, x_n$ of elements of $\mfa$ is called an $\mfa$-filter regular sequence on $M$ if
$$ 
\biggl(\sum_{t=1}^iAx_t \biggr) M :_M x_{i+1}
\subseteq \biggl(\sum_{t=1}^iAx_t \biggr) M :_M \mfa^\infty
$$
for all $0\leq i<n$. When such property holds in any order, we will say that the sequence $x_1,\ldots, x_n$ forms an unconditioned $\mfa$-filter regular sequence on $M$. Note that $x_1,\ldots, x_n$ is a poor $M$-sequence \cite[pp. 39]{sz1} if and only if it is an $A$-filter regular sequence on $M$. Also, note that if $x_1,\ldots,x_n$ is an $\mfa$-filter regular sequence on $M$, then $x_1,\ldots, x_{n-l}$ is a poor $M_{x_n}$-sequence, where $M_{x_n}$ is the ordinary localization of $M$ with respect to $x_n$. 

The following characterization of filter regular sequences is included 
here for the reader's convenience. 

\begin{proposition}{\label{prop11}}
{\cite[2.2]{ns}}\;Let $x_1,\ldots, x_n$ be a sequence of elements of $\mfa$. Then the following statements are equivalent: 

\noindent{\rm (i)} $x_1,\ldots, x_n$ is an $\mfa$-filter regular sequence on M.

\noindent{\rm (ii)} $x_1,\ldots, x_n$ is a poor $M_P$-sequence for all $P\in{\rm Supp}(M)$ with $\mfa\not\subset P$.

\noindent{\rm (iii)} $x_1^{t_1},\ldots, x_n^{t_n}$ is an $\mfa$-filter regular sequence on $M$  for all positive integers $t_1,\ldots, t_n$.
\end{proposition}

As we mentioned in the introduction, to prove the main result of this note, we need to choose a better generating set for the ideal $\mfa$. The following 
proposition indicates the choice.

\begin{proposition}{\label{prop12}}
Suppose that $\mfa$ is generated by $n$ elements. Then $\mfa$ has a sequence of generators of length $n$ which forms an unconditioned $\mfa$-filter regular  sequence on M.
\end{proposition}

\begin{proof} 
Let $\mfa =\overset{n}{\underset{t=1}{\sum}} Aa_t$. Suppose, inductively, that $i$ is an integer with $0\leq i < n$ and it has been proved that there exists an unconditioned $\mfa$-filter regular sequence on $M$ of length $i$, say $x_1,\ldots, x_i$, such that $\mfa =\overset{i}{\underset{t=1}{\sum}} Ax_t+ \overset{n}{\underset{t=i+1}{\sum}} Aa_t$. This is certainly the case when $i=0$. Set 
$$
T = \bigl\{
P \;\bigl|\;  \mfa\not\subset P\text{ and }P\in{\rm Ass}
\hbox{$\bigr(M\bigr/\underset{t\in I}{\sum} x_tM\bigl)$}
\text{ for some }I\subseteq\{1,... ,i\}
\bigr\}.
$$
Then, by \cite[theorem 124]{ka}, there exists $x_{i+1}\in\mfa$ such that $\displaystyle x_{i+1}\not\in\bigcup_{P\in T} P$ and that 
$$
x_{i+1}-a_{i+1}\in\sum_{t=1}^iAx_t+ \sum_{t=i+2}^n Aa_t.
$$ 
Now, it follows from \cite[2.1]{ks} that $x_1,\ldots,x_{i+1}$ is an unconditioned $\mfa$-filter regular sequence on $M$. This completes the inductive step and the result is proved by induction. 
\end{proof}\medskip

Let $x_1,\ldots,x_n$ be a sequence of elements of $A$ and let $\displaystyle\mfa = \sum_{i=1}^n Ax_i$. Write the complexes $C(x,M)^\bullet$. and $C(\mfU(x),M)$  as 
$$
0\overset{f^{-1}}{\lar}C(x,M)^{0}\overset{f^{0}}{\lar}C(x,M)^{1} 
\overset{f^{1}}{\lar}\ldots{\lar} C(x,M)^{i}
\overset{f^{i}}{\lar} C(x,M)^{i+1}{\lar}\ldots
$$
and
$$
0\overset{e^{-1}}{\lar}M \overset{e^{0}}{\lar} U(x)_1^{-1}M
\overset{e^{1}}{\lar}\ldots{\lar}U(x)_i^{-i}M\overset{e^{i}}{\lar} U(x)_{i+1}^{-i-1}M{\lar}
\ldots
$$
respectively. Note that, in view of \cite[3.3 (ii)]{sz1} and \cite[2.1]{sz2}, $U(x)_i^{-i} M = 0$ for all $i > n$. For each $k\in\bbN$, there is an $A$-homomorphism 
$$
\theta_x^k: C(x, M)^k\lar U(x)_k^{-k} M
$$
which is such that, for all $m\in M$ and $\a\in\bbN$,
$$
\theta_x^k\biggl(
\frac{m}{(x_{i(1)},\ldots,x_{i(k)})^\a}
\biggr)
=
\begin{cases}
\frac{m}{(x_1^\a,\ldots,x_k^\a)}&\text{ if }
(i(1),... ,i(k)) = (1,... ,k)
\\ 
0 & \text{ otherwise}.
\end{cases}
$$ 
Let $\theta_x^0$ be the identity map of $M$ and put $\theta_x^i = 0$ for all $i > n$. Then 
$$
\Theta_x = (\theta_x^i)_{i\geq 0}: C(x,M)^\bullet\lar C(\mfU(x),M)
$$
is a homomorphism of complexes. For each $i\geq 0$, let $\theta_x^{i*} :{\rm Ker}\,f^i/{\rm Im}\,f^{i-1} \rar {\rm Ker}\, e^i/{\rm Im}\, e^{i-1}$ denote the homomorphism induced by $\Theta_x$. 

\begin{lemma}{\label{lemma13}}
$\theta_x^{n*}$ is an isomorphism.
\end{lemma}

\begin{proof} 
Clearly $\theta_x^{n*}$ is surjective. Suppose that $a\in M$ and $\a\in\bbN$ are such that 
$$
\theta_x^{n*}\biggl(
{\rm Im}\, f^{n-l} + \frac{a}{(x_1\ldots x_n)^\a}
\biggr)= 0
$$
Then there exist $b\in M$ and $\beta\in\bbN$ such that 
$$ 
\frac{a}{(x_1^\a,\ldots,x_n^\a)}
=
\frac{b}{(x_1^\beta,\ldots,x_{n-1}^\beta,1)}
$$
in $U(x)_n^{-n}M$. Hence, by \cite[1.7]{oc}, there exists $\delta\in\bbN$ with $\delta\geq\max\{\a,\beta\}$ such that
$$ 
x_1^{\delta-\a}\ldots x_n^{\delta-\a}a
=
x_1^{\delta-\beta}\ldots x_{n-1}^{\delta-\beta}x_n^\delta b
+\sum_{i=1}^{n-1}
x_i^\delta m_i
$$
for some $m_1,\ldots,m_{n-1}\in M$. Therefore, in $C(x,M)^n$, we have
$$
\begin{array}{l}\disp
\frac{a}{(x_1\ldots x_n)^\a}
=
\frac{x_1^{\delta-\a}\ldots x_n^{\delta-\a}a}{(x_1\ldots x_n)^\delta}
=
\frac{x_1^{\delta-\beta}\ldots x_{n-1}^{\delta-\beta}x_n^\delta b}{(x_1\ldots x_n)^\delta}+\sum_{i=1}^{n-1}
\frac{x_i^\delta m_i}{(x_1\ldots x_n)^\delta}
\\[1.5ex]\disp 
= 
f^{n-1}\biggl(
\frac{m_1}{(\hat x_1\ldots x_n)^\delta},
\frac{-m_2}{(x_1\hat x_2\ldots x_n)^\delta},
\ldots,\frac{(-1)^{n-2} m_{n-1}}{(x_1\ldots\hat x_{n-1} x_n)^\delta},
\frac{(-1)^{n-1} x_1^{\delta-\beta}\ldots x_{n-1}^{\delta-\beta}b}{(x_1\ldots x_{n-1}\hat x_n)^\delta}
\biggr)
\end{array}
$$
where the character with $\;\hat{}\;$ means that it is deleted, as required. 
\end{proof}\medskip

As we have mentioned in the introduction, the main purpose of this note is to establish a quasi-isomorphism from the complex $C(x, M)^\bullet$ to a certain complex of modules of generalized fractions. Since the homology modules of $C(x,M)^\bullet$ are independent of the choice of sequence of generators for the ideal $\mfa =\overset{n}{\underset{i=1}{\sum}} Ax_i$, we may assume, in view of the proposition \ref{prop12}, that $x_1,\ldots,x_n$ is an $\mfa$-filter regular sequence on $M$. Indeed, under such assumption, we shall prove, by induction, that $\Theta_x$ is a quasi-isomorphism. The next lemma provides a basis for that induction. 

\begin{lemma}{\label{lemma14}}
Let $x_1, x_2$ be an $\mfa = Ax_1 + Ax_2$-filter regular sequence on $M$. Then $\Theta_x$ is a quasi-isomorphism.
\end{lemma}

\begin{proof}
In view of lemma \ref{lemma13}, the only non-trivial point is the proof that $\theta_x^{1*}$ is an isomorphism. To show that $\theta_x^{1*}$  is  injective, let $\theta_x^{1*}\bigl({\rm Im}\,f^0 + (\frac{a_1}{x_1^\a},\frac{a_2}{x_2^\a})\bigr) = 0$. Then, in view of \cite[1.7]{oc}, there exist $\gamma\in\bbN$ and $b\in M$ such that
\begin{align}{\label{eq-1}}
x_1^\gamma(a_1 - x_1^\a b) = 0.
\end{align}
Now, since $x_1,x_2$ is an $\mfa$-filter regular sequence on $M$, it follows from \eqref{eq-1}, in view of proposition \ref{prop11}, that $x_2^\eta(a_1 - x_1^\a b) = 0$  for some $\eta\in\bbN$. Next, since $\bigl(\frac{a_1}{x_1^\a}, \frac{a_2}{x_2^\a}\bigr)\in{\rm Ker}\,f^1$, we have  $(x_1x_2)^\beta(x_1^\a{a_2} - x_2^\a{a_1}) = 0$ for some $\beta\in\bbN$. Now, let $\nu = \max\{\eta,\beta\}$. Then, using the above formulas and proposition  \ref{prop11}, it is easy to see that 
\begin{align}{\label{eq-2}}
x_2^{\tau+\nu}(a_2 - x_2^\a b) = 0
\end{align}
for some $\tau\in\bbN$. Now, it follows from \eqref{eq-1} and \eqref{eq-2} that   $\bigl(\frac{a_1}{x_1^\a},\frac{a_2}{x_2^\a}\bigr)\in{\rm Im}\,f^0$. Therefore ${\rm Ker}(\theta_x^{1*})=0$. 
\end{proof}\medskip

Next we show that $\theta_x^{1*}$ is surjective. To do this, let $\frac{a}{(x_1^\a)}\in{\rm Ker}\, e^1$. Then, by \cite[1.7]{oc}, there exist $\delta\in\bbN$ with $\delta\geq\a$ and $b \in M$ such that  $x_1^{\delta-\a}x_2^\delta{a} = x_1^\delta{b}$. Let $Y= \bigl(\frac{a}{x_1^\a}, \frac{b}{x_2^\delta}\bigr)\in C(x, M)^1$. Then $Y\in  {\rm Ker}\,f^1$ and $\theta_x^{1*}({\rm Im}\,f^0 + Y) = {\rm Im}\,e^0 + \frac{a}{(x_1^\a)}$, as required. 

Next, we establish the main result of this note. In the course of the proof of that result we shall need to use some local cohomology theory. For $i\geq 0$, we use $H^i_\mfb$ to denote the $i$-th right derived functor of the local cohomology functor with respect to an ideal $\mfb$ of $A$. Also, for a complex $X$ of $A$-modules and $A$-homomorphisms, we use $H^i(X)$ to denote the $i$-th homology module of $X$.

\begin{theorem}{\label{thm15}}
Let $\mfa =\overset{n}{\underset{i=1}{\sum}}Ax_i$ and suppose that $x_1,\ldots, x_n$  is an $\mfa$-filter regular sequence on $M$. Then, with the above notations, 
$$\Theta_x: C(x,M)^\bullet\lar C(\mfU(x),M)$$
is a quasi-isomorphism.
\end{theorem}

\begin{proof} 
The proof will appear elsewhere. 
\end{proof}\medskip

\begin{remark}{\label{rmk16}} 
It should be observed that if $M$ is a graded module over the graded ring $A$ and $U$ is a triangular subset of $A^n$ composed of homogeneous elements, then the module $U^{-n}M$ of generalized fractions has graded structure as an $A$-module (see \cite{en}) which is such that, for a homogeneous element $b\in M$ and $(u_1,\ldots, u_n)\in U$, the degree of the fraction $\frac{b}{(u_1,\ldots,u_n)}$ is $\deg(b)-\overset{n}{\underset{i=1}{\sum}}\deg(u_i)$. Thus, for an ideal $\mfa = \overset{n}{\underset{i=1}{\sum}}Ax_i$ of $A$, if $x_1,\ldots,x_n$ is a homogeneous $\mfa$-filter regular sequence on $M$, then the homology modules of the complex $C(\mfU(x), M)$ have graded modules structure. On the other hand, it is standard to use the graded \v Cech complex $C(x, M)^\bullet$ to define a  grading on the local cohomology module $H^i_\mfa(M)$, for $i\geq 0$. However, since the quasi-isomorphism $\Theta_x:C(x,M)^\bullet\rar C(\mfU(x),M)$ is homogeneous, the resulting gradings are always the same. Therefore one can compute the graded local cohomology modules $H^i_\mfa(M)$ via modules of generalized fractions which is perhaps simpler than the use of \v Cech complex. 

The terms of the \v Cech complex $C(x, A)^\bullet$ has the vanishing properties ${\rm Tor}^A_i\bigl(C(x, A)^k, M\bigr) = 0$ for all $i \geq 1$ and all $k$. The following theorem, which describes the behaviour of certain fraction formation along exact sequences, provides similar vanishing properties for the complex $C\bigl(\mfU(x),A\bigr)^\bullet$.
\end{remark}

\begin{theorem}{\label{thm17}}
Let $0 \rar N' \overset{f}{\rar} N\rar N''\rar 0$ be an exact sequence of $A$-modules and $A$-homomorphisms and let $U_n$ be a triangular subset of $A^n$ such that, for each $(x_1,\ldots,x_n)\in U_n$, the sequence $x_1,\ldots,x_{n-1}$ is an $\mfa$-filter regular sequence on $N''$ and $x_n\in\mfa$. Then the sequence
$$
0 \lar U_n^{-n}N' 
\overset{U_n^{-n}f}{\lar} U_n^{-n}N \lar U_n^{-n}N''\lar  0
$$
is exact. 
\end{theorem}

\begin{proof} Let $(x_1,\ldots,x_n)\in U_n$. In view of \cite[3.5 and 2.9]{sz1}, it is enough to show that $U_n^{-n}f$ is injective whenever $U_n=\{(x_1^{\a_1}, \ldots,x_n^{\a_n}) \mid \a_1,\ldots,\a_n\in\bbN\}$. We achieve this by induction on $n(\geq 1)$. The case $n = 1$ is clear. Suppose, inductively, that $n \geq 2$ and the result has been proved for $n-1$. Let $m'\in N'$ and $\a\in\bbN$ be such that 
$$
\frac{f(m')}{(x_1^\a,\ldots,x_n^\a)} = 0
$$
in $U_n^{-n}N$. Then, by \cite[1.7]{oc}, there exists $\delta\in\bbN$ with $\delta\geq\a$ and $m_1,\ldots, m_{n-1}\in N$ such that $x_1^{\delta-\a}\ldots x_n^{\delta-\a} f(m')= x_1^\delta m_1 +\ldots+x_{n-1}^\delta m_{n-1}$. Next, using the above equality in conjuction with the notion of filter regular sequence, one can see that $\disp x_n^\gamma m_{n-1}- f(c) \in \sum_{i=1}^{n-2} x_i^\delta N$ for some $c\in N'$ and $\gamma\in\bbN$. It therefore follows that
$$
x_1^{\delta-\a}\ldots x_{n-1}^{\delta-\a}x_n^{\gamma+\delta-\a}f(m')
- x_{n-1}^\delta f(c) \in\sum_{i=1}^{n-2} x_i^\delta N.
$$
Hence, in view of \cite[3.3 (ii)]{sz1} and the inductive hypothesis, 
$$ 
\frac{x_n^{\gamma+\delta-\a}m'}{(x_1^\a,\ldots,x_{n-1}^\a)}
=
\frac{x_{n-1}^{\delta}c}{(x_1^\delta,\ldots,x_{n-1}^\delta)}
$$
in $U_{n-1}^{-n+1}N'$. Thus $\frac{m'}{(x_1^\a,\ldots,x_{n}^\a)}=0$ by  \cite[2.1]{sz2}, as required. 
\end{proof}\medskip

\begin{consequence}{\label{conseq18}}
Let U be a triangular subset of $A^n$ which consists entirely of $\mfa$-filter regular sequences on $M$ and $A$. Then ${\rm Tor}_i^A(U^{-n}A, M) = 0$ for all  $i\geq 1$.
\end{consequence}

\begin{proof} 
Consider a finite free resolution for M and use the natural equivalence of functors $U^{-n}A\otimes -\cong U^{-n}$ in conjunction with theorem \ref{thm17} to deduce the result.
\end{proof}

\medskip


\begin{thebibliography}{xxx}

\bibitem{bs} Brodmann,~M.\,P., Sharp,~R.\,Y., Local cohomology: {\it an algebraic introduction with geometric applications}, Cambridge University Press, Cambridge, 1998.\smallskip 

\bibitem{en} Enshaei,~R., {\it Modules of generalized fractions and graded rings and modules}, Ph.D. Thesis, University of Sheffield 1987.\smallskip

\bibitem{ka} Kaplansky,~I, {\it Commutative rings}, Allyn\& Bacon: Boston, 1970.\smallskip

\bibitem{ks} Khashyarmanesh,~K., Salarian,~Sh., Zakeri,~H., {\it On the associated graded module of an ideal generated by an unconditioned strong d-sequence}, J.~Math. Kyoto Univ., (1999), {\bf 39} (4), 607-618.\smallskip

\bibitem{ns} Nagel,~U, Schenzel,~P., {\it Cohomological annihilators and Castelnuovo-Mumford regularity}, in {\it Commutative algebra: Syzygies, multiplicities, and birational algebra}, Contemp. Math. Providence, RI 1994, 
307-388.\smallskip

\bibitem{oc} O'Carroll,~L., {\it On the generalized fractions of Sharp and Zakeri}, J.~London Math. Soc., (1983), {\bf 28} (2), 417-427.\smallskip

\bibitem{sz1} Sharp,~R.\,Y., Zakeri,~H., {\it Modules of generalized fractions}, Mathematika, (1982), {\bf 29}, 32-41.\smallskip 

\bibitem{sz2} Sharp,~R.\,Y., Zakeri,~H., {\it Local cohomology and modules of generalized fractions}, Mathematika, (1982), {\bf 29}, 296-306.
\end{thebibliography}
\end{document}